\documentclass[12pt,twoside,leqno]{article}
\usepackage{amsmath}
\usepackage{amssymb}
\usepackage{amsxtra}
\usepackage{amscd}
\usepackage{amsthm}
\usepackage[mathscr]{eucal}
\usepackage{color}
\usepackage{fancyhdr}

\theoremstyle{plain}
\newtheorem{thm}[subsection]{Theorem}
\newtheorem{prop}[subsection]{Proposition}

\newtheorem{lem}[subsection]{Lemma}
\newtheorem{conj}[subsection]{Conjecture}
\newtheorem{ques}[subsection]{Question}%

\theoremstyle{definition}

\newtheorem{rem}[subsection]{Remark}
\newtheorem{para}[subsection]{}

\newtheorem{sbpara}[subsubsection]{}

\newenvironment{pf}{\proof[\proofname]}{\endproof}

\begin{document}

\title{Logarithmic geometry and Frobenius}

\author{Kazuya Kato, Chikara Nakayama, Sampei Usui}
\pagestyle{fancy}
\lhead{}
\rhead{Logarithmic geometry and Frobenius}
\renewcommand{\headrulewidth}{0pt}

\maketitle

\newcommand\Cal{\mathcal}
\newcommand\define{\newcommand}

\define\gp{\mathrm{gp}}%
\define\fs{\mathrm{fs}}%
\define\an{\mathrm{an}}%
\define\mult{\mathrm{mult}}%
\define\Ker{\mathrm{Ker}\,}%
\define\Coker{\mathrm{Coker}\,}%
\define\Hom{\mathrm{Hom}\,}%
\define\Ext{\mathrm{Ext}\,}%
\define\rank{\mathrm{rank}\,}%
\define\gr{\mathrm{gr}}%
\define\cHom{\Cal{Hom}}
\define\cExt{\Cal Ext\,}%

\define\cB{\Cal B}
\define\cC{\Cal C}
\define\cD{\Cal D}
\define\cO{\Cal O}
\define\cS{\Cal S}
\define\cM{\Cal M}
\define\cG{\Cal G}
\define\cH{\Cal H}
\define\cE{\Cal E}
\define\cF{\Cal F}
\define\cN{\Cal N}
\define\fF{\frak F}
\define{\sW}{\Cal W}

\define\Dc{\check{D}}
\define\Ec{\check{E}}

\newcommand{\N}{{\mathbb{N}}}
\newcommand{\Q}{{\mathbb{Q}}}
\newcommand{\Z}{{\mathbb{Z}}}
\newcommand{\R}{{\mathbb{R}}}
\newcommand{\C}{{\mathbb{C}}}
\newcommand{\bN}{{\mathbb{N}}}
\newcommand{\bQ}{{\mathbb{Q}}}
\newcommand{\bF}{{\mathbb{F}}}
\newcommand{\bZ}{{\mathbb{Z}}}
\newcommand{\bP}{{\mathbb{P}}}
\newcommand{\bR}{{\mathbb{R}}}
\newcommand{\bC}{{\mathbb{C}}}
\newcommand{\bG}{{\mathbb{G}}}
\newcommand{\bbQ}{{\bar \mathbb{Q}}}
\newcommand{\ol}[1]{\overline{#1}}
\newcommand{\too}{\longrightarrow}
\newcommand{\respect}{\rightsquigarrow}
\newcommand{\compatible}{\leftrightsquigarrow}
\newcommand{\upc}[1]{\overset {\lower 0.3ex \hbox{${\;}_{\circ}$}}{#1}}
\newcommand{\Gmlog}{\bG_{m, \log}}%
\newcommand{\Gm}{\bG_m}%
\newcommand{\ep}{\varepsilon}
\newcommand{\Spec}{\operatorname{Spec}}
\newcommand{\nilp}{\operatorname{nilp}}
\newcommand{\prim}{\operatorname{prim}}
\newcommand{\val}{{\mathrm{val}}} 
\newcommand{\n}{\operatorname{naive}}
\newcommand{\bs}{\operatorname{\backslash}}
\newcommand{\Gal}{\operatorname{{Gal}}}
\newcommand{\gal}{{\rm {Gal}}({\bar \Q}/{\Q})}
\newcommand{\galp}{{\rm {Gal}}({\bar \Q}_p/{\Q}_p)}
\newcommand{\gall}{{\rm{Gal}}({\bar \Q}_\ell/\Q_\ell)}
\newcommand{\wep}{W({\bar \Q}_p/\Q_p)}
\newcommand{\wel}{W({\bar \Q}_\ell/\Q_\ell)}
\newcommand{\Ad}{{\rm{Ad}}}
\newcommand{\BS}{{\rm {BS}}}
\newcommand{\even}{\operatorname{even}}
\newcommand{\End}{{\rm {End}}}
\newcommand{\odd}{\operatorname{odd}}
\newcommand{\GL}{\operatorname{GL}}
\newcommand{\np}{\text{non-$p$}}
\newcommand{\g}{{\gamma}}
\newcommand{\G}{{\Gamma}}
\newcommand{\Lam}{{\Lambda}}
\newcommand{\La}{{\Lambda}}
\newcommand{\lam}{{\lambda}}
\newcommand{\la}{{\lambda}}
\newcommand{\uL}{{{\hat {L}}^{\rm {ur}}}}
\newcommand{\uQp}{{{\hat \Q}_p}^{\text{ur}}}
\newcommand{\sel}{\operatorname{Sel}}
\newcommand{\dt}{{\rm{Det}}}
\newcommand{\Sig}{\Sigma}
\newcommand{\fil}{{\rm{fil}}}
\newcommand{\SL}{{\rm{SL}}}
\renewcommand{\sl}{{\frak{sl}}}%
\newcommand{\spl}{{\rm{spl}}}
\newcommand{\st}{{\rm{st}}}
\newcommand{\Isom}{{\rm {Isom}}}
\newcommand{\Mor}{{\rm {Mor}}}
\newcommand{\bg}{\bar{g}}
\newcommand{\id}{{\rm {id}}}
\newcommand{\cone}{{\rm {cone}}}
\newcommand{\al}{a}
\newcommand{\ChL}{{\cal{C}}(\La)}
\newcommand{\Image}{{\rm {Image}}}
\newcommand{\toric}{{\operatorname{toric}}}
\newcommand{\torus}{{\operatorname{torus}}}
\newcommand{\Aut}{{\rm {Aut}}}
\newcommand{\Qp}{{\mathbb{Q}}_p}
\newcommand{\barQp}{{\mathbb{Q}}_p}
\newcommand{\Qpur}{{\mathbb{Q}}_p^{\rm {ur}}}
\newcommand{\Zp}{{\mathbb{Z}}_p}
\newcommand{\Zl}{{\mathbb{Z}}_l}
\newcommand{\Ql}{{\mathbb{Q}}_l}
\newcommand{\Qlur}{{\mathbb{Q}}_l^{\rm {ur}}}
\newcommand{\F}{{\mathbb{F}}}
\newcommand{\eps}{{\epsilon}}
\newcommand{\epsLa}{{\epsilon}_{\La}}
\newcommand{\epsLaVxi}{{\epsilon}_{\La}(V, \xi)}
\newcommand{\epsOLaVxi}{{\epsilon}_{0,\La}(V, \xi)}
\newcommand{\Qplin}{{\mathbb{Q}}_p(\mu_{l^{\infty}})}
\newcommand{\otimesQplin}{\otimes_{\Qp}{\mathbb{Q}}_p(\mu_{l^{\infty}})}
\newcommand{\galFl}{{\rm{Gal}}({\bar {\Bbb F}}_\ell/{\Bbb F}_\ell)}
\newcommand{\gallur}{{\rm{Gal}}({\bar \Q}_\ell/\Q_\ell^{\rm {ur}})}
\newcommand{\galFF}{{\rm {Gal}}(F_{\infty}/F)}
\newcommand{\galFv}{{\rm {Gal}}(\bar{F}_v/F_v)}
\newcommand{\galF}{{\rm {Gal}}(\bar{F}/F)}
\newcommand{\epsVxi}{{\epsilon}(V, \xi)}
\newcommand{\epsOVxi}{{\epsilon}_0(V, \xi)}
\newcommand{\plim}{\lim_
{\scriptstyle 
\longleftarrow \atop \scriptstyle n}}
\newcommand{\sig}{{\sigma}}
\newcommand{\ga}{{\gamma}}
\newcommand{\del}{{\delta}}
\newcommand{\Vss}{V^{\rm {ss}}}
\newcommand{\Bst}{B_{\rm {st}}}
\newcommand{\Dpst}{D_{\rm {pst}}}
\newcommand{\Dcrys}{D_{\rm {crys}}}
\newcommand{\DdR}{D_{\rm {dR}}}
\newcommand{\Fin}{F_{\infty}}
\newcommand{\Kla}{K_{\lambda}}
\newcommand{\Ola}{O_{\lambda}}
\newcommand{\Mla}{M_{\lambda}}
\newcommand{\Det}{{\rm{Det}}}
\newcommand{\Sym}{{\rm{Sym}}}
\newcommand{\LaSa}{{\La_{S^*}}}
\newcommand{\cX}{{\cal {X}}}
\newcommand{\MHG}{{\frak {M}}_H(G)}
\newcommand{\tauMla}{\tau(M_{\lambda})}
\newcommand{\Fvur}{{F_v^{\rm {ur}}}}
\newcommand{\Lie}{{\rm {Lie}}}
\newcommand{\cL}{{\cal {L}}}
\newcommand{\cW}{{\cal {W}}}
\newcommand{\fq}{{\frak {q}}}
\newcommand{\cont}{{\rm {cont}}}
\newcommand{\SC}{{SC}}
\newcommand{\Om}{{\Omega}}
\newcommand{\dR}{{\rm {dR}}}
\newcommand{\crys}{{\rm {crys}}}
\newcommand{\hatSig}{{\hat{\Sigma}}}
\newcommand{\rdet}{{{\rm {det}}}}
\newcommand{\ord}{{{\rm {ord}}}}
\newcommand{\BdR}{{B_{\rm {dR}}}}
\newcommand{\BdRO}{{B^0_{\rm {dR}}}}
\newcommand{\Bcrys}{{B_{\rm {crys}}}}
\newcommand{\Qw}{{\mathbb{Q}}_w}
\newcommand{\barkappa}{{\bar{\kappa}}}
\newcommand{\cP}{{\Cal {P}}}
\newcommand{\cZ}{{\Cal {Z}}}
\newcommand{\cA}{{\Cal {A}}}
\newcommand{\oppLa}{{\Lambda^{\circ}}}

\renewcommand{\bar}{\overline}
\newcommand{\et}{\mathrm{\acute{e}t}}
\newcommand{\loget}{\mathrm{log\acute{e}t}}
\newcommand{\pri}{{{\rm prim}}}
\newcommand{\add}{{\rm{add}}}

\begin{abstract} Based on the logarithmic algebraic geometry and the theory of Deligne systems, 
we define an abelian category of $\ell$-adic sheaves with weight filtrations on a logarithmic scheme over a finite field, which is similar to the category of variations of mixed Hodge structure.  
 We consider asymptotic behaviors and simple cases of higher direct images of objects of this category. This category is closely related to the monodromy-weight conjecture.  

\end{abstract}

\renewcommand{\thefootnote}{\fnsymbol{footnote}}

\renewcommand{\thefootnote}{\arabic{footnote}} 

\section*{Introduction}\label{s:intro}
In this paper, based on the theory of logarithmic \'etale topology (\cite{Na0}, 
\cite{Na})  and  the theory of Deligne systems (\cite{SC}, \cite{BPR}), we define an abelian category $\cA_X$ of $\ell$-adic sheaves with weight filtrations on an fs logarithmic scheme $X$ over a finite field.

This category is similar to the category of variations of mixed Hodge structure.  
We describe asymptotic behaviors of objects of $\cA_X$ (Theorem \ref{SL2Th}), which are similar to asymptotic behaviors of variations of mixed Hodge structure in degeneration. We consider simple cases (Propositions \ref{Hi}, \ref{Hi2}) of higher direct images of objects of $\cA_X$.  

This category is closely related to the monodromy-weight conjecture (\ref{conjMW}, \ref{MW}).

In Section \ref{s:1}, we explain our main ideas, main definitions,  and results of this paper. In Section \ref{s:rev}, we consider admissible actions of cones and review the theory of Deligne systems. In Section \ref{s:rat}, we review the space of ratios used in Sections \ref{s:pfab} and \ref{s:SL2}.  In Section \ref{s:pfab}, we discuss basic things about the category $\cA_X$. In Section \ref{s:SL2}, we  state and prove Theorem \ref{SL2Th}. In Section \ref{s:Ex}, we consider higher direct images of objects of $\cA_X$ and study some examples.

  K.\ Kato was 
partially supported by NFS grants DMS 1303421, DMS 1601861, and DMS 2001182.
C.\ Nakayama was 
partially supported by JSPS Grants-in-Aid for Scientific Research (B) 23340008, (C) 16K05093, and (C) 21K03199.
S.\ Usui was 
partially supported by JSPS Grants-in-Aid for Scientific Research (B) 23340008, (C) 17K05200, and  (C) 22K03247.

\section{Logarithmic local systems and Frobenius}\label{s:1}

In this Section \ref{s:1},  we introduce our main ideas and main results. 

Let $k$ be a finite field.

Let $\ell$ be a prime number which is not the characteristic of $k$. We fix an isomorphism of fields $\iota: \bar \Q_{\ell}\cong \C$. 

Let $X$ be an fs log  scheme of finite type over  $k$. 

\begin{para}\label{1.1}

 Let $\cB_X$ be the category of smooth $\bar \Q_{\ell}$-sheaves $H$ on $X$ for the log \'etale topology endowed with a smooth finite increasing  filtration (which we denote by $W$). Let $\cA_X$ be the full subcategory of $\cB_X$ consisting of objects which satisfy the following conditions 1 and 2 for every closed point $x$ of $X$.

1. The local monodromy of $H$ at $x$ is unipotent and is admissible with respect to $W$. Here we use the fixed isomorphism $\iota : \bar \Q_{\ell}\cong \C$ to define the admissibility (\ref{adm}) of the action of $\R_{\geq 0}$-cones.

2. For the monodromy cone $\sig(x)$ at $x$, for the relative monodromy filtration $W(\sig(x))$ and  for every $w\in \Z$, $\gr^{W(\sig(x))}_wH$ is of Frobenius weight $w$ as a representation of $\Gal(\bar k/k)$. Here we use the fixed isomorphism $\iota$ to define the weight of the Frobenius. (For the monodromy cone $\sig(x)$, see \ref{lp}. For $W(\sig(x))$, see \ref{adm}.)
\end{para}

\begin{rem}\label{dep} 
(1) Not only the Frobenius weight in the condition 2, the admissibility in the condition 1 depends on the choice of the isomorphism $\iota$ (see \ref{depend}). If $X$ is of log rank $\leq 1$, this admissibility is independent of the choice. 

(2) There are two kinds of log \'etale topologies: the kummer log \'etale topology considered in \cite{Na0} and the full log \'etale topology considered in \cite{Na}.  
  We can use any one of them, that is, 
the theory in this paper becomes the same for both 
($\cB_X$ are the same, $\cA_X$ are the same, and so on) by \cite{Na} Theorem 5.17. 
  But to fix ideas, we use the kummer log \'etale topology below.
\end{rem} 

\begin{para}\label{abelian}

 We hope that $\cA_X$ is an analogue of the category of variations of mixed Hodge structure and has similar properties to the latter. 
 
For example, 
the category $\cA_X$ is an abelian category. This follows from the theory of Deligne (\cite{BPR} Lemma 6.33) on Deligne systems, as is explained in \ref{pfab}. 

 We have a theorem (Theorem \ref{SL2Th}) on the asymptotic behavior of an object of $\cA_X$, which  is similar to the asymptotic behavior of a variation of mixed Hodge structure in degeneration. This is obtained from the theory of Deligne systems by using the space of ratios in Section \ref{s:rat}.  Theorem \ref{SL2Th} is useful in the study of the asymptotic behaviors of the  regulator maps and height pairings in degeneration. This study is illustrated in \cite{Ka} Remark 2.4.18 and will be explained elsewhere.

The category $\cA_X$ is related to the monodromy-weight conjecture as in the following Conjecture \ref{conjMW} and \ref{MW}. 

\end{para}

\begin{conj}\label{conjMW} Let $f: X\to Y$ be a projective log smooth vertical saturated morphism of fs log schemes over $k$. Then  $R^mf_*(\bar \Q_{\ell})$ is an object of $\cA_Y$ of pure $W$-weight $m$.

\end{conj}

\begin{para}

By proper base change theorem (\cite{Na0} Theorem (5.1), \cite{Na} Theorem 6.1), Conjecture \ref{conjMW} is reduced to the case where $Y$ is an fs log point.
\end{para}

\begin{para}\label{MW}
 In the case where $Y$ is the standard log point, Conjecture \ref{conjMW} implies the monodromy weight conjecture. In fact, if $K$ is a non-archimedean local field with residue field $k$ and if $\frak X$ is a projective scheme over $O_K$ with semi-stable reduction, then $R^mf_*\bar \Q_{\ell}$ for $X=\frak X\otimes_{O_K} k$ and $Y=\Spec(k)$ with the canonical log structures is identified with the representation $H^m_{\et}({\frak X}\otimes_{O_K} {\bar K}, \bar \Q_{\ell})$ of $\Gal(\bar K/K)$, and the monodromy-weight  conjecture for $\frak X$ is exactly that $R^mf_*\bar \Q_{\ell}$ is an object of $\cA_Y$ and is pure of $W$-weight $m$. 
 
 \end{para}
 
\begin{rem}
(1) Conjecture \ref{conjMW} becomes not true if we replace the assumption of the projectivity by the properness. In fact, for the special fiber  $X\to Y$ of a formal model of the  Hopf surface over a non-archimedean local field ($X$ is algebraic though the Hopf surface is not algebraic but only analytic), $R^1f_*\bar \Q_{\ell}$ is $\bar \Q_{\ell}$ of weight $0$, not of weight $1$. 

(2) Without the condition saturated, we may not have the unipotence of local monodromy.
\end{rem}

\begin{para}\label{lp} We give some pictures of the categories $\cB_X$ and  $\cA_X$ in the case where $x$ is a point.

Let $X$ be an fs log point whose underlying scheme is $\Spec(k)$ ($k$ is a finite field). We  denote $X$ also by the small letter $x$ because it is a point. Let $\bar k$ be a separable closure of $x$ and let $\bar x=\Spec(\bar k)$. We have an exact sequence 
$$\begin{CD}
1 @>>> \Gal(\bar x(\log)/\bar x) @>>> \Gal(\bar x(\log)/x) @>>> \Gal(\bar x/x) @>>> 1.\\
@. @| @| @| \\
@. \Hom((M^{\gp}_X/\cO_X^\times)_{\bar x},\Z)\otimes \hat \Z(1)'   @. \pi_1^{\log}(x)  @. \Gal(\bar k/k)
\end{CD}
$$
Here $\bar x(\log)$ is $\Spec(\bar k)$ endowed with the log structure obtained by adding $n$-th roots of the log structure of $\bar x$ for all $n\geq 1$ which are not divisible by the characteristic of $k$ (i.e., $\bar x(\log)$ is the associated log separably closed field in the sense of \cite{Na0} Definition (2.5)), and 
$\hat \Z(1)'$ is the product of $\Z_{\ell'}(1)$, where $\ell'$  ranges over all prime numbers which are different from the characteristic of $k$. If $q$ denotes the order of the finite field $k$ and if $F$ is an element of $\pi_1^{\log}(x)$ whose image in $\Gal(\bar k/k)$ is the $q$-th power map, we have $F\gamma F^{-1}=\gamma^q$ for all $\gamma\in \Gal(\bar x(\log)/\bar x)$. 

By taking the stalk at $\bar x(\log)$, an object of  $\cB_X$ is identified with a finite dimensional representation of $\pi_1^{\log}(x)$ over $\bar \Q_{\ell}$ endowed with a $\pi_1^{\log}(x)$-stable finite increasing filtration $W$. The action of $\Gal(\bar x(\log)/\bar x)$ on this stalk is called the local monodromy at $x$. 

The cone $\sig(x):= \Hom((M_X/\cO_X^\times)_{\bar x},\R^{\add}_{\geq 0})$, where $\R_{\geq 0}^{\add}=\{x\in\R\;|\; x \geq 0\}$ with the additive structure,  is called the {\em monodromy cone} at $x$. The admissibility of the local monodromy in \ref{1.1} means  the admissibility of the logarithm action of this monodromy cone.

\end{para}

\begin{para}\label{stlp1} 

Assume that $X=\{x\}$ is the standard log point $\Spec(k)$ whose log structure is associated to $\N\to k\;;\;n\mapsto 0^n$.  

Then $\Gal(\bar x(\log)/x)\cong \Gal(K_{\mathrm{tame}}/K)$ for a non-archimedean local field $K$ with the residue field $k$ and for its  maximal tame extension $K_{\mathrm{tame}}$.

Let $N\in \sig(x)$ be the standard generator of
$\Hom((M_X/\cO_X^\times)_{\bar x},\bN)=\Hom(\bN,\bN)=\bN$. 

For $r\geq 0$, we define an object $S_r$ of $\cB_X$. Fix a splitting of the surjective homomorphism $\pi_1^{\log}(x) \to \Gal(\bar k/k)$. 
Let 
\begin{equation*}\tag{1}
S_1=\bar\Q_{\ell} \oplus \bar\Q_{\ell}(-1)
\end{equation*}
\noindent as a representation of $\Gal(\bar k/k)$. Let $e_1=(1,0)\in S_1$ and let $e_2=(0, e)$, where $e$ is any fixed base of $\bar\Q_{\ell}(-1)$. 
  Define the action of $\pi_1^{\log}(x)$ on $S_1$ by $N(e_2)=e_1$ and $N(e_1)=0$. Define $W$ of $S_1$ to be pure of weight $1$. Then $S_1$ is an object of $\cA_X$. 
 Let $S_r$ be the $r$-th symmetric power of $S_1$. This is an object of $\cA_X$ of pure $W$-weight $r$.
 
 As Proposition  \ref{stlp2} below shows,  $S_r$ are simple objects of $\cA_X$. This is remarkable because we have for example  an exact sequence
 $$0 \to \bar \Q_{\ell} \to S_1 \to \bar \Q_{\ell}(-1)\to 0$$
 of $\ell$-adic sheaves on the log \'etale site of  $X$ which tells that $S_1$ is not simple as an $\ell$-adic sheaf.

For $w\in \Z$, let $\cA_w$ be the full subcategory of $\cA_X$ consisting of all objects which have pure  $W$-weight $w$.

On the other hand, for a smooth $\bar \Q_{\ell}$-sheaf $H$ on $X$ with no log structure, if $H$ has pure Frobenius weight $w$,  we regard $H$ as an object of $\cA_X$ with pure $W$-weight $w$ in the natural way. 
Let $\cA_w'$ be the category of families $(H_r)_{r\geq 0}$ of smooth $\bar \Q_{\ell}$-sheaves  $H_r$ on $X$ with no log structure such that $H_r$ has pure Frobenius weight $w-r$ and such that $H_r=0$ for almost all $r$.  We have a functor
\begin{equation*}\tag{2}%
\Phi: \cA_w' \to \cA_w\;;\; (H_r)_r \mapsto \bigoplus_r S_r\otimes H_r.
\end{equation*}
\end{para}

\begin{prop}\label{stlp2}  Assume that $X$ is the standard log point.

$(1)$  The functor $\Phi$ gives an equivalence of categories 
$$\cA_w' \overset{\simeq}\to \cA_w.$$

$(2)$ Let $(H_r)_r$ be an object of $\cA'_w$ and let $H=\Phi((H_r)_r)$. Then $H$ is semi-simple in $\cA_X$ if and only if $H_r$ are semi-simple for all $w$. $H$ is simple in $\cA_X$ if and only if  there is an $r$ such that $H_r$ is of rank $1$ and such that $H_s=0$ for all $s\neq r$. 

$(3)$ A $W$-pure object $H$ of $\cA_X$ is semi-simple if and only if $H$ is semi-simple as a representation of $\Gal(\bar k/k)$ (forgetting the local monodromy). 
\end{prop}

This will be proved in \ref{stlp3}.

\begin{para}\label{RandQ}  We are using the $\R_{\geq 0}$-cone $\sig(x)= \Hom((M_X/\cO_X^\times)_{\bar x},\R^{\add}_{\geq 0})$. If we use the $\Q_{\geq 0}$-cone $\sig_{\Q}(x):= \Hom((M_X/\cO_X^\times)_{\bar x},\Q^{\add}_{\geq 0})$ instead of $\sig(x)$, we have a similar but different theory. Concerning this, see \ref{RQ}.  
\end{para}

\begin{para} The category $\cA_X$ has a crystalline analogue and a  $p$-adic Hodge analogue. We will consider them in a sequel of this paper. 

\end{para}

\begin{rem} One of the authors (K. Kato) published a paper \cite{Ka:KJM} concerning Deligne systems. It is shown in \cite{BPR} Section 6 that this paper \cite{Ka:KJM} is wrong.  He hopes that  \cite{Ka:KJM} is never used by any person. He hoped that the paper \cite{Ka:KJM} could be used in the study of degeneration of a family of motives over a non-archimedean local field (\cite{Ka:KJM} 1.9). We hope that Theorem \ref{SL2Th} in our present paper is on the right way for that purpose.

\end{rem}

\section{Reviews on admissible actions of cones and on Deligne systems}\label{s:rev}

We review relative monodromy filtration, admissibility of an action of a cone, Deligne splitting, and Deligne systems (See \cite{SC}).

\begin{para}\label{2.4.b} 
Let  $V$ be a  vector space  
endowed with a finite increasing filtration $W$, and let $N: V \to V$ be a nilpotent linear map such that $NW_w\subset W_w$ for any $w\in \Z$. A finite  increasing filtration $\sW$ on $V$  is called {\em the relative monodromy filtration of $N$ with respect to $W$} if it satisfies the following two conditions (i) and (ii).

(i) $N\sW_w\subset \sW_{w-2}$ for any $w\in \Z$.

(ii) For every $w\in \Z$ and every integer $r\geq 0$, 
the map  $N^r: \gr^{\sW}_{w+r}\gr^W_w \overset{\cong}\to \gr^{\sW}_{w-r}\gr^W_w$ is an isomorphism.

The relative monodromy filtration $\sW$ of $N$ with respect to $W$ is unique if it exists (\cite{De} 1.6.13).

\end{para}

\begin{para}\label{sigma} By a {\em finitely generated cone}, we mean a subset $\sigma$ of a finite dimensional  $\R$-vector space $V$ such that  $\sigma=\R_{\geq 0}N_1+\dots+ \R_{\geq 0}N_n$ for some $N_1,\dots, N_n\in V$. Here $\R_{\geq 0}=\{x\in \R\;|\;x\geq 0\}$. 

We denote the $\R$-linear span of $\sig$ in $V$ by $\sig_\R$. 
\end{para}

\begin{para}\label{adm} Let  $\sigma$ be a finitely generated sharp cone and  
let $V$ be a finite dimensional $\C$-vector space endowed with a finite increasing filtration $W$. By an {\em admissible action} of $\sigma$ on $V$, we mean an $\R$-linear map $h: \sigma_\R\to \End_\C(V)$ satisfying the following conditions (i)--(iii).

(i) $h(N)h(N')=h(N')h(N)$ for all $N,N'\in \sigma_\R$ and $N(W_w)\subset W_w$ for all $N\in \sig_{\bR}$ and $w\in \Z$.

(ii) For every $N\in \sigma_\R$, $h(N)$ is nilpotent.

(iii) There is a family $(W(\tau))_{\tau}$ of finite increasing filtrations $W(\tau)$ on $V$, where $\tau$ ranges over all faces of $\sig$,  satisfying the following conditions (iii-1), (iii-2), (iii-3).

(iii-1)  For every face $\tau$ of $\sigma$ and for every $N\in \sigma_\R$ and $w\in \Z$,  we have $N(W(\tau)_w)\subset W(\tau)_w$.

(iii-2) $W(\{0\})=W$.

(iii-3) For any faces $\tau$, $\tau'$, $\tau''$ of $\sigma$  such that $\tau\supset \tau'$ and  for any $w\in \Z$ and any element $N$ of the interior of $\tau$, 
the restriction of $W(\tau)$ to $W(\tau'')_w$  is the relative monodromy filtration of the restriction of $N$ to $W(\tau'')_w$ with respect to the restriction of $W(\tau')$ to $W(\tau'')_w$. 

\end{para}

\begin{rem}
Concerning the above condition (iii):

(1) By the conditions (iii-2) and (iii-3), for any face $\tau$ of $\sigma$ and for any element $N$ of the interior of $\tau$, $W(\tau)$ is the relative monodromy filtration of $N$ with respect to $W$. Hence the family $(W(\tau))_{\tau}$ is unique if it exists. 

(2) (1) tells that for any face $\tau$ of $\sigma$ and for any $N\in \tau_{\R}$ and $w\in \Z$,  we have $N(W(\tau)_w)\subset W(\tau)_{w-2}$. 
  This is because the interior of $\tau$ generates $\tau_{\R}$ as an $\R$-linear space.

\end{rem}

\begin{para}\label{act}  For an admissible action $\sig$ on $V$ and for a face $\tau$ of $\sig$, we have $NW(\tau)_w\subset W(\tau)_w$ for all $w\in \Z$ and for all $N\in \sig_\R$. 

In fact,  for an element $N'$ of the interior of $\tau$, $W(\tau)$ is the relative monodromy filtration of $N'$ with respect to $W$, and hence $\exp(N)W(\tau)$ is the relative monodromy filtration of $\exp(N)N'\exp(-N)=N'$ with respect to $\exp(N)W=W$, and hence coincides with $W(\tau)$. 

\end{para}

\begin{para} The admissibility is stable under taking direct sums, tensor products, and duals. 

\end{para}

  \begin{para}\label{Nd} We review the theory of Deligne splitting (\cite{SC} Theorem 1) due to Deligne. 
  
  We return to the setting in \ref{2.4.b}. Assume that 
the relative monodromy  filtration $\sW$ of $N$ with respect to $W$ exists, and assume that we are given a splitting $Y$  of $\sW$ satisfying the following two conditions. 

(i) $N$ is of weight $-2$ for $Y$.

(ii) $Y$ is compatible with $W$. That is, the action of $\bG_m$ on $V$ associated to $Y$ ($a\in \bG_m$ acts on the part of $V$ of weight $w$ for $Y$ by $a^w$) keeps $W$. 

Then there is a unique splitting $Y'$ of $W$, which we call the {\em Deligne splitting}, satisfying the following two conditions. 
\medskip

(i) $Y'$ is compatible with $Y$. That is, there is  an action of $\bG_m^2$ on $V$ in which the first $\bG_m$ gives $Y'$ and the second $\bG_m$ gives $Y$.

(ii) Let 
$$\End(V)\cong \bigoplus_{w,m} \;\gr^{\sW}_m\gr^W_w\End(V)$$
be the isomorphism given by $(Y',Y)$, and  write the image of $N\in \End(V)$ in $\bigoplus_{w,m} \;\gr^{\sW}_m\gr^W_w\End(V)$ as $\sum_{d\geq 0} N_d$ with $N_d\in \gr^{\sW}_{-2}\gr^W_{-d}\End(V)$. Then for $d\geq 1$, $N_d$ belongs to the primitive part of 
 $\gr^{\sW}_{-2}\gr^W_{-d}\End(V)$, that is, $N_d$ is killed by $\Ad(N)^{(d-1)/2}: \gr^{\sW}_{-2}\gr^W_{-d}\End(V) \to \gr_{-2d}^{\sW}\gr^W_{-d}\End(V)$.

\end{para}

\begin{para}\label{D1} In \ref{D1} and \ref{D2}, we review the theory of Deligne systems (\cite{SC} Section 2, \cite{BPR} Section 6).

An {\em $n$-variable  Deligne system} is a system $(V, (W^j)_{0\leq j\leq n}, (N_j)_{1\leq j\leq n},  Y)$,  
where $V$ is a finite dimensional vector space,  $W^j$  ($0\leq j\leq n$) are finite increasing filtrations on $V$, $N_j: V \to V$ ($1\leq j\leq n$) are mutually commuting nilpotent linear operators, and $Y$ is a splitting of  $W^n$  satisfying the following conditions (i)--(iii). 

(i) For $1\leq j \leq n$, $0\leq i \leq n$ and $w\in \Z$, the restriction of $W^j$ to $W^i_w$ is  the relative monodromy filtration of the restriction of $N_j$ to $W^i_w$ with respect to the restriction of $W^{j-1}$ to $W^i_w$. 

(ii)  $N_i\in W^j_0\End(V)$ for all $i$ and $j$, and 
$N_i\in W^j_{-2}\End(V)$ if $i \leq j$. 

(iii) $N_j$ are of weight $-2$ for $Y$  for all $1\leq j\leq n$, and $Y$ is compatible with $W^j$ for all $0\leq j\leq n$.

\medskip

The category of $n$-variable Deligne systems is an abelian category (\cite{BPR} Lemma 6.33, due to Deligne).

\end{para}

\begin{para}\label{D2} (\cite{SC} Theorem 2, due to Deligne.)

Let $(V, (W^j)_{0\leq j\leq n}, (N_j)_{1\leq j\leq n}, Y)$ be an $n$-variable Deligne system. Then there is a unique action $\rho$ of 
$\bG_m\times \SL(2)^n$ on $V$ characterized by the properties (i) and (ii) below.

For $0 \leq j \leq n$, define the splitting $Y^j$ of $W^j$ by downward induction on $j$ in the following way. $Y^n:=Y$. For $1\leq j<n$, since  $W^{j+1}$ is the relative monodromy filtration of $N^{j+1}$ associated to $W^j$, the splitting $Y^{j+1}$ of $W^{j+1}$ (given by downward induction on $j$) and $N_{j+1}$ define a splitting of $W^{j}$ by the theory of Deligne splitting (\ref{Nd}). For each $1\leq j\leq n$, let $\hat N_j$ be the component of $N_j$ of degree $0$ for the splittings $Y^i$ of $W^i$ for $1\leq i<j$.

{\rm (i)} The Lie action $\Lie(\rho)$ sends the matrix $\begin{pmatrix} 0 & 1 \\ 0&0\end{pmatrix}$ in the $j$-th $\sl(2)$ in $\sl(2)^n$ to $\hat N_j$.

{\rm (ii)} For $0\leq j\leq n$, let $\tau_j$ be the action of $\bG_m$ on $V$  given by the splitting $Y^j$ of $W^{j}$. Then $\bG_m$ in $\bG_m\times \SL(2)^n$ acts via $\tau_0$, and for $1\leq j\leq n$, 
$\tau_j(a)=\tau_0(a)\rho(b)$ for $a\in \bG_m$.
  Here $b$ is the element of $\SL(2)^n$ whose $i$-th component is $\mathrm{diag}\,(a^{-1},a)$ if $1\leq i\leq j$ and is $1$ if $j<i\leq n$.

 \end{para}

\begin{para}\label{adD} Let $\sig$ and $V$ be as in \ref{adm} and assume that we are given an admissible action of $\sig$ on $V$. Assume that we are given a sequence of faces 
$\sig_0\subsetneq \sig_1\subsetneq \dots \subsetneq \sig_n=\sig$, an element $N_j$ of the interior of $\sig_j$ for each $1\leq j\leq n$, and a splitting $Y$ of $W(\sigma)$  such that $N_j$ are of weight $-2$ for $Y$ for $1\leq j\leq n$ and $Y$ is compatible with $W(\sig_j)$ for $0
\leq j\leq n$.

Then we have a Deligne system $(V, (W^j)_{0\leq j \leq n}, (N_j)_{1\leq j\leq n},Y)$, where $W^j:=W(\sig_j)$, and hence an action of $\bG_m \times \SL(2)^n$ on $V$ by \ref{D2}. 

\end{para}

\begin{para}\label{tau} In \ref{adD}, let $\tau$ be a face of $\sig$. Then the action of $\bG_m \times \SL(2)^n$ on $V$ keeps $W(\tau)$.

In fact, for each $w\in \Z$ and for $U:=W(\tau)_w$, $(U, (W^j|_U), (N_j|_U), Y|_U)$ is a Deligne system by the condition (iii-3) in \ref{adm}, and the associated action of $\bG_m \times \SL(2)^n$ on $U$ is compatible with the action of $\bG_m \times \SL(2)^n$ on $V$ by the characterizations of these actions.

\end{para}

\section{Reviews on the space of ratios}\label{s:rat}

We review the space of ratios defined in \cite{KNU4} Section 4 and used in the study of degenerating Hodge structures.

\begin{para} Definition (\cite{KNU4} 4.1.3). Let  $\cS$ be a sharp fs monoid. We denote the semi-group law of $\cS$ multiplicatively. 

The space of ratios $R(\cS)$ is the set of all maps $(\cS\times \cS)\smallsetminus \{(1,1)\}\to[0, \infty]$ satisfying the following conditions (i)--(iii).

\medskip

(i) $r(f,g)= r(g,f)^{-1}$.

(ii) $r(f,g)r(g,h)=r(f,h)$ if $\{r(f,g), r(g,h)\}\neq \{0, \infty\}$. 

(iii) $r(fg,h)=r(f, h)+r(g,h)$.

\medskip

We endow $R(\cS)$ with the topology of simple convergence. It is compact.

\end{para}

\begin{para}\label{rank2} Example. $R(\N^2)$ is homeomorphic to the interval $[0, \infty]$ as a topological space. In fact, if $(q_j)_{j=1,2}$ denotes the standard base of $\N^2$, $r\in R(\N^2)$ corresponds to $r(q_1, q_2)\in [0, \infty]$. 

\end{para}

\begin{para}\label{1:1} There is a bijection between $R(\cS)$ and the set of all equivalence classes of families $((\sig_j)_{1\le j \le n}, (N_j)_{1\le j \le n})$,
where $n\geq 0$, $\sig_j$ are faces of the cone $\sig:= \Hom(\cS, \R_{\geq 0}^{\add})$
such that $\{0\}:=\sig_0\subsetneq \sig_1\subsetneq \dots \subsetneq \sig_n=\sig$ and $N_j$ is an element of the interior of $\sig_j$.  Two such families 
$((\sig_j)_j, (N_j)_j)$ and $((\sig'_j)_j, (N'_j)_j)$ are equivalent if and only if $\sig'_j=\sig_j$ for $1\leq j\leq n$ and there are $c_j\in \R_{>0}$ such that $N_j'\equiv c_jN_j\bmod \sig_{j-1,\R}$ for $1\leq j\leq n$.   (Cf.\ \cite{KNU4} 4.1.6.) 

In fact, for such a family $((\sig_j)_j, (N_j)_j)$, the corresponding $r\in R(\cS)$ is as follows.  Let $\cS=\cS^{(0)}\supsetneq \cS^{(1)}\supsetneq \dots \supsetneq  \cS^{(n)}=\{1\}$ be the sequence of faces of $\cS$ for which $\sig_j$ is the annihilator of $\cS^{(j)}$. Then for $(f,g)\in \cS\times \cS\smallsetminus \{(1,1)\}$, if $j$ is the smallest integer such that $\{f,g\}$ is not contained in 
$\cS^{(j)}$,  $r(f,g)= N_j(f)/N_j(g)$.

\end{para}

\begin{para}\label{n=1} Let $R(\cS)_1$ be the subset of $R(\cS)$ consisting of all elements such that the corresponding sequences of faces of $\sig$ (\ref{1:1}) are of length one, that is, $n=1$. 

If $\cS=\{1\}$,  $R(\cS)_1$ is empty.

 Assume $\cS\neq \{1\}$, and let $\sig^{\circ}$ be the interior of $\sig$. Then, we have a bijection $\psi:\sig^{\circ}/\R_{>0} \to R(\cS)_1$ defined by $\psi(N)(f,g)= N(f)/N(g)$ $(N \in \sig^{\circ})$, and $\psi(N)$ corresponds in \ref{1:1} to the family $(\sig, N)$. 

If $\cS\neq \{1\}$, $R(\cS)_1$ is a dense open subset of $R(\cS)$. The density is proved as follows. Let $r\in R(\cS)$ be the class of $((\sig_j)_{1\le j \le n}, (N_j)_{1\le j \le n})$. Then $r$ is the limit of $\psi(\sum_{j=1}^n y_jN_j)$, where $y_j\in \R_{>0}$ and $y_j/y_{j+1}\to \infty$ for $1\leq j<n$. The openness is shown as follows. Assume that $f_1, \dots, f_m$ ($f_j \in \cS\smallsetminus \{1\}$) generate $\cS$. Then $r\in R(\cS)$ belongs to $R(\cS)_1$ if and only if $r(f_i, f_j)\in (0, \infty)$ for every $i,j$.

\end{para}

\begin{para}\label{x:}  If $x$ is an fs log point, we denote the space $R((M_x/\cO_x^\times)_{\bar x})$ by $x_{[:]}$ and  $R((M_x/\cO_x^\times)_{\bar x})_1$ by $x_{[:],1}$.

 In the case where the underlying scheme of $x$ is $\Spec(\C)$, 
the space
$x_{[:]}$ 
 is important in log Hodge theory to treat $\SL(2)$-orbits.  See \cite{KNU4} Theorem 4.5.2 and \cite{KNU6} Section 4. 

\end{para}

\begin{para}\label{val}  Let $x$ be an fs log point. To work on the topological space $x_{[:]}$ as we will do in Section \ref{s:SL2} may give the impression that we work in the category of topological spaces leaving algebraic geometry. But this is not a correct feeling. If $\bar x_{\val}$ denotes the inverse limit of log blowing ups of $\bar x$, we have a surjective continuous map $\bar x_{\val}\to x_{[:]}$ such that the topology of $x_{[:]}$ is the quotient topology of the topology of $\bar x_{\val}$ (\cite{KNU4} 4.1.11, 4.1.12). Hence for an object $H$ of $\cA_x$, the behavior of $H$ on the topological space $x_{[:]}$ gives information of the behavior of $H$ on the space ${\bar x}_{\val}$ of algebraic nature. Since $\bar x$ and $\bar x_{\val}$ are identical for the log \'etale topology (the full log \'etale topology (\cite{Na}), not the kummer log \'etale topology), the behavior of $H$ on the topological space $x_{[:]}$ gives information on the behavior of $H$ on the log \'etale site of $\bar x$.

\end{para}

\section{The category $\cA_X$}\label{s:pfab}

We discuss basic facts about the category $\cA_X$.     
  
  \begin{lem}\label{act2} Let $x$ be an fs log scheme whose underlying scheme is $\Spec(k)$ for a finite field $k$. Let $H$ be an object of $\cA_x$ and  let $V$ be the stalk of $H$.

 Let $\tau$ be a face of $\sig(x)$, and let $W(\tau)$ be the weight filtration on $V$ associated to $\tau$ ($\ref{adm}$). Then $W(\tau)$ is  invariant under the action of $\pi_1^{\log}(x)$.

  \end{lem} 
  
  This is a stronger version  of \ref{act} in the present situation. 
  
  \begin{pf} Let $q$ be the order of the finite field $k$, and let $F\in \pi_1^{\log}(x)$ be an element whose image in $\Gal(\bar k/k)$ is the $q$-th power map.  Let $N$ be an element of the interior of $\tau$. Then since $W(\tau)$ is the relative monodromy filtration of $N$ with respect to $W$, $FW(\tau)$ is the relative monodromy filtration of $FNF^{-1}$ with respect to $FW=W$. Since $FNF^{-1}=qN$ (\ref{lp}), $FW(\tau)$ coincides with $W(\tau)$. 
  
  Since these $F$ generate the group $\pi_1^{\log}(x)$, $W(\tau)$ is invariant under the action of $\pi_1^{\log}(x)$.
\end{pf}
    
    \begin{lem}\label{Frobspl} Let $x$, $k$, $H$ and $V$ be as in Lemma $\ref{act2}$. 
     Fix a splitting of $\pi_1^{\log}(x)\to \Gal(\bar k/k)$, and let $Y$ be the splitting of $W(\sig(x))$ by the Frobenius weights for the action of $\Gal(\bar k/k)$.  
  Then the following hold. 
  
  $(1)$ For every $N\in \sig(x)_\R$, $N:V\to V$ is of weight $-2$ for $Y$.

  $(2)$ For every face $\tau$ of $\sig(x)$,  $Y$ is compatible with $W(\tau)$.
    
    \end{lem}
    
    \begin{pf} Let $F$ be the $q$-th power Frobenius in $\Gal(\bar k/k)$ which acts on $V$ via the splitting.
    
    (1) follows from $FNF^{-1}=qN$.

    We prove (2). For each $c\in \C^\times$, if $V_c$ denotes the generalized eigenspace of $F$ in $V$ for the eigenvalue $c$, the projection $V\to V_c$ is given by a polynomial in $F$. By Lemma \ref{act2}, this projection sends $W(\tau)_w$ into $W(\tau)_w$ for each $w$, and hence $W(\tau)_w= \bigoplus_c\; W(\tau)_w\cap V_c$. Hence $Y$ is compatible with $W(\tau)$. 
    \end{pf}
  
  \begin{para}\label{SL2A} Let the situation be as in Lemma \ref{Frobspl}. 
    
  Assume that we are given a sequence of faces $\sig_0\subsetneq \sig_1\subsetneq \dots \subsetneq \sig_n=\sig$ and an element $N_j$ of the interior of $\sig_j$ for each $1\leq j\leq n$. 
  Here $\sig_0$ need not be $\{0\}$.

Then by \ref{adD} and Lemma \ref{Frobspl}, we have a Deligne system $(V, (W^j)_{0\leq j \leq n}, (N_j)_{1\leq j\leq n},Y)$, where $W^j:=W(\sig_j)$, and hence an action of $\bG_m \times \SL(2)^n$ on $V$ by \ref{D2}.  By \ref{tau}, for every face $\tau$ of $\sig(x)$ and for every $w\in \Z$, $W(\tau)_w$ is stable under this action of $\bG_m \times \SL(2)^n$.

  \end{para}

\begin{para}\label{SL2Del} Let the situation be as in Lemma \ref{Frobspl}. We consider the space of ratios $x_{[:]}$. Let $\mu\in x_{[:]}$.

  Recall that $\mu$ is the equivalence class of a family 
$((\sig_j)_{1 \leq j \leq n}, (N_j)_{1\leq j\leq n})$, where $\sig_j$ are faces of $\sig(x)$ such that $\{0\}=\sig_0\subsetneq \sig_1\subsetneq \dots \subsetneq \sig_n=\sig(x)$ and $N_j$ is an element of the interior of $\sig_j$ (see \ref{1:1}).  

As in \ref{SL2A}, this family 
$((\sig_j)_{1 \leq j \leq n}, (N_j)_{1\leq j\leq n})$ 
determines a representation of $\bG_m \times \SL(2)^n$.

By the characterization of Deligne splitting reviewed in \ref{Nd}, for $1 \leq j \leq n$, the splitting $Y^{j-1}$ of $W^{j-1}$  in \ref{D2} depends only on $N_j \bmod \sig_{j-1,\R}$.  Because elements of $\sig_{j-1,\R}$ have only weights $\leq -2$ parts for $Y^{j-1}$,  $\hat N_j$ in \ref{D2} depends only on
$N_j\bmod \sig_{j-1,\R}$. Hence this action of $\bG_m \times \SL(2)^n$  is determined by $((\sig_j)_{1\leq j\leq n}, (N_j\bmod \sig_{j-1,\R})_{1\leq j\leq n})$.

If we replace $N_j$ by $a_jN_j$ ($1\leq j \leq n)$ for $a_j\in \R_{>0}$, $\rho$ is changed by its conjugate by 
 the action of the element of $\SL(2,\C)^n$ whose $j$-th component is $\text{diag}(\sqrt{a_j}, 1/\sqrt{a_j})$. Thus $\mu$ determines the representation of $\bG_m \times \SL(2)^n$ modulo this conjugacy.
 
 This is very similar to the relation between the space of ratios and $\SL(2)$-orbits in Hodge theory  in \cite{KNU4} Theorem 6.3.1 (1), \cite{KNU6} Section 4.

 \end{para}

\begin{para}\label{pfab}
Let $X$ be an fs log scheme over a finite field.

We prove that $\cA_X$ is an abelian category. This is closely related to the fact that the category of $n$-variable Deligne systems is an abelian category (\ref{D1}) and proved in a similar way. 
We also prove the following  (1) and (2).

\medskip

(1) If $H \to H'$ is a morphism in $\cA_X$, then at each closed point $x\in X$,
for every face $\sig$ of the monodromy cone  $\sig(x)$
and $w\in \Z$, the image of $W(\sig)_wH$ in $H'$ coincides with the intersection of the image of $H$ and $W(\sig)_wH'$.

(2) For any exact sequence $0\to H' \to H \to H'' \to 0$ in $\cA_X$, at each closed point $x\in X$,  the sequence $0\to W(\sig)_wH'\to W(\sig)_wH \to W(\sig)_wH'' \to 0$ is exact for every $w$ and for every  $\sig$ as in (1). 

\medskip 
In the following proof, when we work at a closed point $x$ of $X$, we will denote the stalks of $H$, $H'$, etc., just by the same notation $H$, $H'$, etc.

 We prove (1). Apply \ref{SL2A} to the case where $n=1$, $\sig_0=\sig$, $\sig_1=\sig(x)$, and $N_1$ is an element of the interior of $\sig(x)$. Then this defines actions of $\bG_m$ on $H$ and $H'$ splitting $W(\sig)$, and the homomorphism $H\to H'$ is compatible with these actions. 
This proves (1).

 We prove that $\cA_X$ is an abelian category. 

First, for a morphism $H\to H'$ of $\cA_X$, we prove that the kernel $K$ with the induced $W$ and the cokernel $C$ with the induced $W$ are objects of $\cA_X$. We consider the kernel. For $x\in X$ and for each face $\tau$ of $\sig(x)$, let $W(\tau)K$ be the filtration on $K$ induced by the filtration $W(\tau)H$. We prove that we have an admissible action with the family
 $(W(\tau)K)_{\tau}$ of 
relative monodromy filtrations.  It is sufficient to prove that the condition (iii-3) in \ref{adm} is satisfied. 
Let $\tau, \tau', \tau''$ be as in this condition (iii-3).

Apply \ref{SL2A} to the case where $n=2$, $\sig_0=\tau'$, $\sig_1=\tau$, $\sig_2=\sig(x)$, $N_1$ is an element of the interior of $\tau$, and $N_2$ is an element of the interior of $\sig(x)$.
Then we have  actions of $\SL(2)^2$ on $H$ and $H'$, and the homomorphism $H\to H'$ is compatible with these actions. This induces an action of $\SL(2)^2$ on $K$ and  on $W(\tau'')_wK$ for $w\in \Z$ (\ref{SL2A}). Hence for every $r, s\in \Z$ with $r\geq 0$, the map $N_1^r={\hat N}_1^r: \gr^{W(\tau)}_{s+r} \gr^{W(\tau')}_s W(\tau'')_wK \to \gr^{W(\tau)}_{s-r} \gr^{W(\tau')}_s W(\tau'')_wK $ is an isomorphism.

The proof for the cokernel $C$ is similar.

Next, we prove that the map from the image to the coimage  is an isomorphism. But this follows from (1) for $W=W(\{0\})$. This completes the proof of the statement that $\cA_X$ is an abelian category.

 We prove (2). Apply \ref{SL2A} to the case where $\sig_0=\sig$, $\sig_1=\sig(x)$, and $N_1$ is an element of the interior of $\sig(x)$. Then we have 
  the actions of $\bG_m$  on $H$, $H'$, $H''$ splitting $W(\sig)$ and we have  the exact sequence of representations of $\bG_m$. 
Hence the induced sequence of $W(\sig)_w$ is exact.

\end{para}

\begin{para}\label{stlp3}  We prove Proposition \ref{stlp2}.

We prove (1). We write here the proof for the case $w=0$. The method of the proof for the general case is the same, but we assume this just to make the notation simple. We fix a splitting of $\pi_1^{\log}(x)\to \Gal(\bar k/k)$. Hence for every object $H$ of $\cA_X$, we have a splitting of the Frobenius weight filtration. For an object $H$ of $\cA_X$, let $H=\bigoplus_w H^{[w]}$ be the decomposition, where $H^{[w]}$ is the part of $H$ of Frobenius weight $w$ (it is a sub $\bar \Q_{\ell}$-sheaf of $H$). We give the converse functor $\cA_0 \to \cA'_0$ in two ways. First, define $\Psi^-: \cA_0\to \cA'_0$ as $\Psi^-(H)= (H_r)_{r\geq 0}$, where $H_r=\text{Ker}(N: H^{[-r]}\to H^{[-r-2]})$. Next, let $\Psi^+(H)= (H_r)_{r\geq 0}$, where $H_r= \text{Ker}(N^{r+1}: H^{[r]}\to H^{[-r-2]})(r)$ (this is the so-called  primitive part). Then we have a natural isomorphism $\Psi^+\overset{\cong}\to \Psi^-$ given by $N^r: \text{Ker}(N^{r+1}: H^{[r]}\to H^{[-r-2]})(r) \overset{\cong}\to \text{Ker}(N: H^{[-r]}\to H^{[-r-2]})$. The composition $\Psi^- \circ \Phi$ is identified with the identity functor of $\cA'_0$. We show that the composition $\Phi \circ \Psi^+$ is isomorphic to the identity functor of $\cA_0$. For an object $H$ of $\cA_X$, if $(H_r)_{r\geq 0}= \Psi^+(H)$, we have the natural isomorphism $\bigoplus_r S_r\otimes H_r\overset{\cong}\to H$  given by
$S_r\otimes H_r\to H$ which sends $j!^{-1}e_1^ie_2^j\otimes x$ ($i+j=r$) to $N^i(x)$. Here we identify $\bar \Q_{\ell}(i)$ with $\bar \Q_{\ell}$ by $e^{\otimes -i} \mapsto 1$ for all $i\in \Z$, where $e$ is the base of $\bar \Q_{\ell}(-1)$ which we fixed in \ref{stlp1}.

(2) follows from (1). 

 (3) follows from (1) and (2). 
\end{para}

\begin{para}\label{depend} We show that the admissibility in the condition 1 in \ref{1.1} depends on the choice of the isomorphism $\iota$.

Let $X$ be an fs log point $(\Spec(k), \N^2 \oplus \cO_X^{\times})$. 
  Define an object $H$ of $\cB_X$ as follows. As a representation of $\Gal(\bar k/k)$, $H=\bar{\Q}_{\ell}\oplus \bar{\Q}_{\ell}(-1)$ with the base $e_1,e_2$ as in \ref{stlp1}. $W$ of $H$ is pure of weight $1$. Take non-zero elements $a_j$ ($j=1,2$) of $\bar{\Q}_{\ell}$. Let the monodromy operators $N_j$ ($j=1,2$) be $N_j(e_1)=0$ and $N_j(e_2)=a_je_1$. Then the local monodromy is admissible if and only if $x_1\iota(a_1) +x_2\iota(a_2)\neq 0$ for every $x_1,x_2\in \R_{>0}$. This condition is equivalent to the condition that $\iota(a_1a_2^{-1})$ is not a negative real number, and hence depends on the choice of $\iota$ (for example, it depends if $a_1a_2^{-1}$ is a square root of $2$).

In this situation,  the local monodromy is admissible for every choice of $\iota$ if and only if $a_1a_2^{-1}$ is a totally positive algebraic number. 
\end{para}

\begin{para}\label{RQ} If we use the $\Q_{\geq 0}$-cone $\sig_{\Q}(x)$ (\ref{RandQ}) instead of the $\R_{\geq 0}$-cone $\sig(x)$, and consider the admissibility using only element of $\sig_{\Q}(x)$, we do not need to use an isomorphism $\bar{\Q}_{\ell}\cong \C$ to define the admissibility. We still can prove that the version of $\cA_X$ for this formulation is an abelian category (the proof is essentially identical with the proof given above). However, this weaker admissibility does not give Theorem \ref{SL2Th}. For example, if $x$ is an fs log point with the log structure $\N^2$ as in \ref{depend}, this weak admissibility cannot give the behavior of   $H$ around the point of $x_{[:]}\cong [0, \infty]$ (\ref{rank2}) corresponding to $\sqrt{2}\in [0, \infty]$. 
This is not nice for the applications explained in \ref{nonarch} below. 
This is the reason why we like to use the $\R_{\geq 0}$-cone.

\end{para}

\section{Asymptotic behaviors}\label{s:SL2}

We prove a non-archimedean analogue Theorem \ref{SL2Th} of the asymptotic behaviors of mixed Hodge structures in degeneration. The latter is related to the SL(2)-orbit theorem in Hodge theory. It has been studied by many people, for example, as in
\cite{CKS},   \cite{Kas}, \cite{KNU0}, \cite{KNU2}, \cite{Sch},  etc.
  
Theorem \ref{SL2Th} is applied to non-archimedean geometry as in  \ref{nonarch} below.

\begin{para}\label{t2}

Let $x$ be an fs log scheme whose underlying scheme is $\Spec(k)$ for a finite field $k$. We assume that the log structure of $x$ is not trivial.

Let $H$ be an object of $\cA_x$ and let $V$ be the stalk of $H$.

We use  the space $x_{[:]}$ of ratios and its dense open subset $x_{[:],1}$ (\ref{x:}). Let $\mu\in x_{[:]}$. We consider the behavior of $H$ when $\nu\in x_{[:],1}$ converges to $\mu$.

\end{para}

\begin{para}\label{Ntnu}

 Let $$E=\End(V)$$ be the set of all linear operators on $V$ regarded as a $\C$-algebra.

Assume  that $\mu\in x_{[:]}$ is the class of $((\sig_j)_{1\leq j\leq n}, (N_j)_{1\leq j\leq n})$  (\ref{1:1}).

Let $M_{\bar x}= M_{\bar x}^{(0)}\supsetneq M_{\bar x}^{(1)}\supsetneq \dots \supsetneq M^{(n)}_{\bar x}= \cO^\times_{\bar x}=k^\times$ be the sequence of faces of $M_{\bar x}$ which is dual to the filtration $(\sig_j)_j$ on $\sig(x)$. Fix $f_j\in M_{\bar x}$ ($1\leq j\leq n$) which is in $M_{\bar x}^{(j-1)}$  but not in $M_{\bar x}^{(j)}$.

Replacing $N_j$ ($1\leq j\leq n$) by its $\R_{>0}$-multiple, we assume that $N_j:(M_x^{\gp}/\cO_x^\times)_{\bar x}\to \R$ sends $f_j$ to $1$ for $1\leq j\leq n$. Fix a splitting of $\pi_1^{\log}(x) \to \Gal(\bar k/k)$. As in \ref{SL2Del}, we have a splitting $Y^j$ of $W^j=W(\sig_j)$ for $0\leq j\leq n$ and we have an action of $\bG_m\times \SL(2)^n$ on $V$.  Let $\hat N_j\in E$ ($1\leq j\leq n$) be as in \ref{D2}. 
For $1\leq j\leq n$, let $\tau_j$ be the action of $\bG_m$ on $V$ corresponding to $Y^j$. 
For $\nu\in x_{[:], 1}$, we define  $$N(\nu)\in E, \quad t(\nu)\in E^\times$$ as follows.

Let $N_{\nu}$ to be the unique element of the interior of $\sig(x)$ such that $\nu$ is the class of $(\sig(x), N_{\nu})$ and such that  $N_{\nu}:(M_x^{\gp}/\cO_x^\times)_{\bar x}\to \R$ sends $f_n$ to $1$. Let $N(\nu)$ be the image of $N_{\nu}$ in $E$.

 Define $$t(\nu):= \prod_{j=1}^{n-1} \tau_j(\nu(f_{j+1}, f_j)^{1/2})\in E^\times.$$

\end{para}

\begin{thm}\label{SL2Th} 
When $\nu\in x_{[:], 1}$  converges to $\mu\in x_{[:]}$, $t(\nu)^{-1}N(\nu)t(\nu)$ converges to $\sum_{j=1}^n \hat N_j$ in $E$.

\end{thm}

This is an analogue of \cite{KNU2} Theorem 2.4.2 (ii) on SL(2)-orbits of mixed Hodge structures.

 The following lemma is actually a part of the theory of Deligne systems, but we present it here with proof because it plays a key role below. 
  
  \begin{lem}\label{Nwt-2}   Let $1\leq j\leq n$ and let $N$ be an element of $\sig_{j,\R}$. Then $N$ is purely of weight $-2$ for $Y^i$ if $j\leq i\leq n$. 
\end{lem}

\begin{pf} 
 This is because $N$ is a homomorphism of Deligne systems 
$(V, (W^i)_{j\leq i\leq n}, (N_i)_{j<i\leq n}, Y^n) \to (V, (W^i)_{j\leq i\leq n}, (N_i)_{j<i\leq n}, Y^n)(-1)$, where $(-1)$ denotes the Tate twist, 
 and the construction of the splittings $Y^i$ $(j \le i \le n)$ is functorial.
 \end{pf}

\begin{para} We prove Theorem \ref{SL2Th}.

  Let $f_{j,\la} \in M^{(j-1)}_{\bar x}$ ($1 \le j \le n$) be the elements such that for any $j$, $(f_{j,\la})_\la$ is a $\bQ$-basis of 
$\bQ \otimes (M^{(j-1),\gp}_{\bar x}/M^{(j),\gp}_{\bar x})$. 
  Let $(N_{j,\la})_{j,\la}$ be the dual base of $(f_{j,\la}\bmod \cO_{\bar x}^\times)_{j,\la}$ in $\sig(x)_\R$. 
  Then $$N(\nu)=\sum_{j,\la}\; \nu(f_{j,\la}, f_n) N_{j,\la}.$$

It is enough to prove that for each $1\leq j\leq n$,  $\Ad(t(\nu))^{-1}\sum_{\la}\;\nu(f_{j,\la}, h) N_{j,\la}$ converges to $\hat N_j$ in $E$. 
  
  We have 
  
  (1)  $N_j\equiv \sum_{\la}\; \mu(f_{j,\la}, f_j) N_{j,\la} \bmod \sig_{j-1,\R}$. 
  
   Write $N_{j,\la}=N^{(0)}_{j,\la}+ N^{(<0)}_{j,\la}$, where $N_{j,\la}^{(0)}$ is of $\tau_i$-weight $0$ for all $i$ such that $1\leq i<j$, and $N_{j,\la}^{(<0)}$ is of $\tau_i$-weight $\leq 0$ for all $i$ such that $1\leq i<j$ and of $\tau_i$-weight $<0$ for some $i$ such that $1\leq i<j$.  By (1), we have
   
   (2) $\hat N_j= \sum_{\la} \mu(f_{j,\la}, f_j) N^{(0)}_{j,\la}$.

When $\nu$ converges to $\mu$,   $\nu(f_{j,\la}, f_j)$ converges to $\mu(f_{j,\la}, f_j)$.

Write  $t(\nu)=t_{\geq j}(\nu)t_{<j}(\nu)$  with $$t_{\geq j}(\nu):= \prod_{i=j}^{n-1} \tau_i(\nu(f_{i+1}, f_i)^{1/2}), \quad t_{<j}(\nu):=\prod_{i=1}^{j-1}\tau_i(\nu(f_{i+1},f_i)^{1/2}).$$

By Lemma \ref{Nwt-2}, $\tau_i(a)^{-1}N_{j,\la}= a^2N_{j,\la}$ if $j\leq i\leq n$. Hence
\begin{align*}
&\Ad(t_{\geq j}(\nu))^{-1} \sum_{\la} \nu(f_{j,\la}, f_n)N_{j,\la}\\
= &\Bigl(\prod_{i=j}^{n-1} \nu(f_{i+1}, f_i)\Bigr)\sum_{\la} \nu(f_{j,\la},f_n)N_{j,\la} =\sum_{\la} \nu(f_{j,\la}, f_j)N_{j,\la} \\
=& \sum_{\la} \nu(f_{j,\la}, f_j)N_{j,\la}^{(0)}+ \sum_{\la} \nu(f_{j,\la}, f_j)N^{(<0)}_{j,\la}. 
\end{align*}
We have $\Ad(t_{<j}(\nu))^{-1}N_{j,\la}^{(0)}=N_{j,\la}^{(0)}$, 
 and $\Ad(t_{<j}(\nu))^{-1}N_{j,\la}^{(<0)}$ tends to $0$ when $\nu$ converges to $\mu$. Hence $\Ad(t(\nu))^{-1}\sum_{\la} \nu(f_{j,\la},f_n)N_{j,\la}$ converges  to $\hat N_j$.

\end{para}

\begin{para}\label{nonarch} Here,  by using an example, we describe how Theorem \ref{SL2Th} is applied to the non-archimedean geometry.

Let $K$ be a non-archimedean local field with finite residue field and let $\pi$ be a prime element of $K$.

We consider the following example.  Endow $\Spec(O_K[t])$ with the log structure generated by $\pi$ and $t$, and assume that our $x$ is the closed point of $\Spec(O_K[t])$ at which $\pi$ and $t$ have value  $0$. 
Let $\frak X$ be an open subscheme of $\Spec(O_K[t])$ containing $x$ and let $U$ be the inverse image of $\frak X$ in $\Spec(K[t, t^{-1}])$. Assume that $H$ comes from a smooth $\bar \Q_{\ell}$-sheaf $\tilde H$ on $\frak X$. Then the restriction of $\tilde H$ to $U$ is a smooth $\bar \Q_{\ell}$-sheaf on the usual \'etale site. 

If  $\alpha\in \bar K^\times$ is near to $0\in \bar K$, $\alpha$ defines a closed point $\Spec(K(\alpha))$ of  $U$ which we denote by $\alpha$. Let $\tilde H(\alpha)$ be the pullback of $\tilde H$ to $\alpha$. We are interested in how the monodromy operator of $\tilde H(\alpha)$ behaves when $\alpha$ tends to $0$. 

Let $x(\alpha)$ be the closed point of $\Spec(O_{K(\alpha)})$ with the standard log structure. If $\alpha$ is near to $0$, the ring homomorphism $O_K[t]\to O_{K(\alpha)}\;;\ t\mapsto \alpha$ induces a morphism $x(\alpha)\to x$ of log schemes. The monodromy cone $\sig(x(\alpha))$ is of rank one, and the homomorphism $\sig(x(\alpha))\to \sig(x)$ sends a non-zero element to an element of the interior of $\sig(x)$. This determines a point $\nu(\alpha)$ of $x_{[:], 1}$ (\ref{n=1}, \ref{x:}). 

In  the homeomorphism $x_{[:]}\cong [0, \infty]\;;\; r \mapsto r(t, \pi)$ (\ref{rank2}), $\nu(\alpha)$ is identified with $v(\alpha)\in [0, \infty]$, where $v$ is the valuation of $K(\alpha)$ normalized by $v(\pi)=1$. 

Since $\tilde H(\alpha)$ is identified with  the pullback of $H$ to $x(\alpha)$ as a representation of $\Gal(K(\alpha)_{\text{tame}}/K(\alpha))=\Gal(\bar \alpha(\log)/\alpha)$, the stalk of $\tilde H(\alpha)$ and the stalk of $H$ are identified. We apply Theorem \ref{SL2Th} by taking $\mu$ to be the point of $x_{[:]}$ corresponding to $\infty\in [0, \infty]$. Then when $\alpha$ tends to $0$, $\nu(\alpha)$ tends to $\mu$. Take $f_1=t$, $f_2=\pi$ in \ref{Ntnu}. Then $N(\nu(\alpha))$ in \ref{Ntnu} is identified 
with the action of the unique element $N_{\alpha}$ of $\sig(x(\alpha))$ which sends the class of $\pi$ in  $M_{x(\alpha)}/\cO^\times_{x(\alpha)}$ to $1$. 

Hence Theorem \ref{SL2Th} describes the behavior of the  monodromy operator $N_{\alpha}$ of $\tilde H(\alpha)$ when $\alpha$ tends to $0$.

The asymptotic behaviors given  by Theorem \ref{SL2Th} are useful for the study of the asymptotic behaviors of the non-archimedean components of regulators and height pairings in degeneration. This is because the non-archimedean components of the regulators and height pairings are described by monodromy operators.

\end{para}

\begin{rem} (1) Slightly refining the above proof of Theorem \ref{SL2Th}, we can show that the function $\nu\mapsto t(\nu)^{-1}N(\nu)t(\nu)$ on $x_{[:],1}$ extends to a complex analytic function in $\nu(f_{j,\la}, f_j)^{1/2}$ and $\nu(f_{j+1}, f_j)^{1/2}$ on an open neighborhood of $\mu$ in $x_{[:]}$. We do not pursue this here. 

(2) This Section \ref{s:SL2} could be presented as a story of the fs monoid $(M_x/\cO^\times_x)_{\bar x}$ and an admissible representation over $\C$ of its dual cone, and of its space of ratios, forgetting the fs log scheme $x$ and  the prime number $\ell$. But the authors prefer the above presentation, for they hope to apply the result to $\ell$-adic sheaves which come from motives over an algebraic variety over a non-archimedean local field, by the method described in \ref{nonarch}.

\end{rem}

\section{Higher direct images}\label{s:Ex}

\begin{para}\label{MWplan} For the Riemann hypothesis part of the Weil conjecture, the proofs  given in \cite{De} and in \cite{Fa} are to consider higher direct images of  mixed sheaves and reduce the problem to the study of higher direct images in the case of relative dimension one. We hope the study of higher direct images for the category $\cA_X$ is important and that the monodromy-weight conjecture is reduced to the case of relative dimension one. 

Here we present an attempt in this direction. 
\end{para}

\begin{para}\label{As} Let $X$ be an fs log scheme of finite type over a finite field.

Let $\cA_X^{s}$ be the full subcategory of $\cA_X$ consisting of objects $H$ such that $\gr^W_wH$ are semi-simple for all $m$.

\end{para}

If the statement in the following Question \ref{Hope} is true, our plan of the study in \ref{MWplan} would be a 
smooth route.

\begin{ques}\label{Hope} Is the following statement true? 

Let $X\to Y$ be a vertical projective log smooth saturated morphism of fs log schemes of finite type over a finite field.  Then for every $i\geq 0$,  the $i$-th direct image functor $H^i$ sends $\cA_X^s$ to $\cA_Y^s$. Here we define $W_wH^i(H)$ as the image of $H^i(W_{w-i}H)$. It sends pure objects of weight $w$ to pure objects of weight $w+m$. 
\end{ques}

In the rest of this Section \ref{s:Ex}, we consider the case where $Y$ is a standard log point and $X$ is a degenerate elliptic curve (\ref{A5}) hoping that it will be helpful for the future studies of Question \ref{Hope}. 
  Concerning this $X \to Y$, we give Proposition \ref{Hi} treating a case for which the statement in Question \ref{Hope} is true, and give a related result Proposition \ref{Hi2}. We give also an example of $H$ in \ref{A4}  which does not belong to  $\cA^s_X$ 
and 
 whose $H^1$ does not belong to $\cA_Y$. 

\begin{para}\label{A5} We consider the following $X\to Y$.  Let $K$ be a non-archimedean local field with finite residue field $k$, let $Y=\Spec(k)$ with the standard log structure, and let $X$ be the special fiber of a projective model over $O_K$ with semi-stable reduction of a Tate elliptic curve over $K$. Regard $X$ as a log smooth vertical saturated fs log scheme over $Y$ in the canonical way.

\end{para}

\begin{prop}\label{Hi} Let $X\to Y$ be as in $\ref{A5}$. 
Let $H$ be an object of $\cA_X$ such that for every $w\in \Z$, $\gr^W_wH$ is the pullback of an object of $\cA_Y$. 

$(1)$ The higher direct images of $H$ belong to $\cA_Y$. 

$(2)$ If  $H$  belongs to $\cA^s_X$, the higher direct images belong to $\cA^s_Y$. 

\end{prop}

We expect that the proof of Proposition \ref{Hi} given below works 
to extend Proposition \ref{Hi} to a general $X\to Y$ of relative dimension one.  See Remark  \ref{rel1}. 

\begin{para}\label{loget} 

Let $X\to Y$ be as in \ref{A5}. We denote $Y$ also by $y$.

In \ref{loget} and \ref{pi10}, we give preparations on the log \'etale cohomology (\cite{Na}) and the log fundamental group (\cite{Ka2} Section 10), respectively.

Let 
$$T^*:=H^1_{\loget}(X \times_Y \bar y(\log), \Z_{\ell}),$$ which is a free $\Z_{\ell}$-module of rank $2$ with an action of $\pi^{\log}_1(y)=\Gal(\bar y(\log)/y)$. 
Let $$T:=\Hom_{\Z_{\ell}}(T^*, \Z_{\ell}), \quad L:= T \otimes_{\Z_{\ell}} \bar \Q_{\ell}, \quad L^*:= T^*\otimes_{\Z_{\ell}} \bar \Q_{\ell}.$$ 
We denote the group law of $T$ multiplicatively, but the group law of $T^*$, $L$, $L^*$ additively.
We regard $L$ as the Lie algebra  over $\bar \Q_{\ell}$ of  the $\ell$-adic Lie group $T$. It is a commutative Lie algebra, and the inclusion map $T\to L$ is thought as the logarithm.

We have a $\Z_{\ell}$-base $(\gamma_j)_{j=1,2}$ of $T$, a lifting $F\in \pi_1^{\log}(y)=\Gal(\bar y(\log)/y)$ of the $q$-th power map in $\Gal(\bar k/k)$, and a topological generator $\gamma_0$ of $\Gal(\bar y(\log)/\bar y)$ such that the action of $\pi_1^{\log}(y)$ on $T$ satisfies
$$F(\gamma_1)=\gamma_1, \quad F(\gamma_2)= \gamma_2^q, \quad \gamma_0(\gamma_1)= \gamma_1\gamma_2, \quad \gamma_0(\gamma_2)=\gamma_2,$$
where $q$ is the order of the finite field $k$.
As representations of $\pi_1^{\log}(y)$,  we have an isomorphism $L^* \cong S_1$, where $S_1$ is as in  \ref{stlp1}, 
such that the base $(e_j)_j$ of $S_1$ in \ref{stlp1} is the image of the dual base of $(\gamma_j)_{j=1,2}$ in $T^*$ by this isomorphism.  

\end{para}

\begin{para}\label{pi10} Let $X\to Y$ be as above. 

We consider the logarithmic  fundamental group $\pi_1^{\log}(X)$. Fixing a closed point $x$ of $X$ and considering the stalk  at $\bar x(\log)$ lying over $\bar y(\log)$, a
smooth $\bar \Q_{\ell}$-sheaf on the log \'etale site $X_{\loget}$ of $X$ is identified with a representation of $\pi_1^{\log}(X)$ over $\bar \Q_{\ell}$. 

We have an exact sequence 

\begin{sbpara}\label{pi1} \;\; $1 \to \pi_1^{\log}(X\times_Y \bar y(\log)) \to \pi_1^{\log}(X) \to \pi_1^{\log}(y)\to 1$, \end{sbpara}

\noindent
which splits ($\pi_1^{\log}(X)$ is the semi-direct product). We have a surjective homomorphism $\pi_1^{\log}(X\times_Y \bar y(\log))\to T$ whose kernel is a pro-$\ell'$ group, where a pro-$\ell'$ group means the inverse limit of finite groups whose orders are coprime to $\ell$.
 Thus we have a surjective homomorphism from $\pi_1^{\log}(X)$ to the semi-direct product of $T$ and $\pi_1^{\log}(y)$ whose kernel is a pro-$\ell'$ group. 
 
\end{para}

\begin{para}\label{howto}

We  discuss how to compute the higher direct images. 
 We use the Lie algebra cohomology.
When we consider a unipotent representation of $T$, we consider the Lie action of the Lie algebra $L$ by taking the logarithm of the action of $T$.

Assume that the action of $\pi_1^{\log}(X \times_Y \bar y(\log))$ on (the stalk of) $H$ is unipotent. 
Let

\begin{sbpara}\label{cpx} \;\; $0 \to H \to H\otimes_{\bar \Q_{\ell}} L^* \to H \otimes_{\bar \Q_{\ell}} \wedge^2 L^* \to 0$\end{sbpara}

\noindent
be the standard complex to compute the Lie algebra cohomology of $L$ with coefficients in $H$. Here we identify $H$ with its stalk endowed with the action of $\pi_1^{\log}(X)$. 
The first $H$ is put in degree $0$ of the complex. This \ref{cpx} is a complex of representations of $\pi_1^{\log}(X)$ (which acts on $L^*$ through $\pi_1^{\log}(X)\to \pi_1^{\log}(y)$) and hence is a complex of $\ell$-adic sheaves on $X_{\loget}$.

As an $\ell$-adic sheaf,  the $i$-th higher direct image $H^i(H)$ of $H$ on $Y$ is identified, by taking the pullback from $Y$  to $X$, with  the $i$-th cohomology of this complex.

\end{para}

\begin{para}\label{Hipf} We prove Proposition \ref{Hi}.  

We identify $L^*$ with the object $S_1$ of $\cA_Y$ pure of $W$-weight $1$
 (\ref{loget})  and regard it as an object of $\cA_X$ by pullback. Then each term of the complex \ref{cpx} with the $W$ of the tensor product is regarded as an object of $\cA_X$. Furthermore, by the assumption 
on $\gr^WH$ in Proposition \ref{Hi}, morphisms in \ref{cpx} are morphisms in $\cA_X$. Thus \ref{cpx} is a complex in $\cA_X$.

As an object of $\cB_Y$, the $i$-th higher direct images of $H$ on $Y$ is identified  with  the $i$-th cohomology of this complex.  
Since $\cA_X$ is an abelian category,  these higher direct images are objects of $\cA_X$. We prove that this object on $Y$ belongs to $\cA_Y$. Take  a closed point  $x$ of $X$ whose log structure  is strict over $Y$. 
Since the local monodromy of the higher direct image at $x$ is unipotent, the local monodromy of the higher direct images on $Y$ is unipotent. Since   the relative monodromy filtration of a non-trivial element of $\sig(x)$ on the higher direct image (regarded as an object of $\cA_X$) is  the Frobenius weight filtration and since $\sig(x) \overset{\cong}\to \sig(y)$, 
the relative monodromy filtration of a non-trivial element of $\sig(y)$ on the higher direct image is the Frobenius weight filtration. This shows that the higher direct images belong to $\cA_Y$.

We prove (2) of Proposition \ref{Hi}. We may  assume that $H$ is $W$-pure and hence $H$ comes from a semi-simple object of $A_Y$ by pullback. Since $L^*\cong S_1$ and $S_r\otimes S_1\cong S_{r+1}\oplus S_{r-1}$ for $r\geq 0$,  Proposition \ref{stlp2} shows that $H\otimes\wedge^i L^*$ are semi-simple for all $i$. Since the $i$-th higher direct image of $H$ is a subquotient of $H\otimes \wedge^i L^*$, it is semi-simple.
\end{para}

\begin{rem}\label{rel1} It seems that the above proof of Proposition \ref{Hi} works for the general case of relative dimension one. We take $L$ to be the Lie algebra of an $\ell$-adic nilpotent quotient of $\pi_1^{\log}(X \times_Y \bar y(\log))$ (then $L$ need not be commutative in this general situation, and $L^*$ becomes of $W$-weight $\geq 1$, not necessarily of $W$-weight $1$). However, for the proof of this generalization, it seems that we have to discuss the $\cA_X$-version of the theory of mixed Hodge structures on  nilpotent quotients of the fundamental groups of an algebraic varieties in \cite{HZ}.   We hope to discuss this elsewhere. 
\end{rem}

An interesting point of the following proposition is that the projectivity of $X\to Y$ which appears in Question \ref{Hope} also appears in the condition (iii) in (3).

\begin{prop}\label{Hi2} Let $X\to Y$ be as in $\ref{A5}$ and let
$H$ be a $W$-pure simple object of $\cA_X$. 

$(1)$ There is $c\in \bar \Q_{\ell}^\times$ such that on the stalk of $H$, the action of $\gamma_1-c$ is nilpotent. 

$(2)$ If $c$ in $(1)$ is not $1$, the higher direct images $H^i(H)$ vanish for all $i$. 

$(3)$ If $c$ in $(1)$ is $1$, the following four conditions are equivalent. 

{\rm (i)} The higher direct images $H^i(H)$ belong to $\cA_Y$ for all $i$.

{\rm (ii)} $H^0(H)$ belongs to $\cA_Y$.

{\rm (iii)} The Lefschetz class map $H^0(H) \to H^2(H(1))$ is an isomorphism of $\ell$-adic sheaves on $Y$. 

{\rm (iv)} $H$ is the pullback of an object of $\cA_Y$.

If these equivalent conditions are satisfied, {\rm (i)} is satisfied with $\cA_Y$ replaced by $\cA^s_Y$ and the map in {\rm (iii)}  becomes an isomorphism in $\cA_Y$.

\end{prop}

We give preparatory lemmas to prove Proposition \ref{Hi2}. 

\begin{lem}
\label{gamma2} Let $X\to Y$ be as in $\ref{A5}$ and let $H$ be an object of $\cA_X$. Then the action of $\gamma_2$ on the stalk of $H$ is unipotent. 

\end{lem}

\begin{pf} Let $\alpha$ be a singular point of $X$. We show that $\gamma_2$ comes from the local monodromy group at $\alpha$. Because the local monodromy of $H$ is unipotent, this will show that the action of $\gamma_2$ is unipotent. We prove that $\gamma_2$ is in the image of  $\pi_1^{\log}(\alpha\times _y \bar y(\log)) \to T$. This map is dual to  $T^*\to H^1_{\loget}(\alpha\times_y \bar y(\log),\Z_{\ell}) \cong \Z_{\ell}$, which sends $e_2$ to $1$ and $e_1$ to $0$. Hence the last map is induced by $\gamma_2$. 
\end{pf}

\begin{lem}\label{eigen}

Let $X\to Y$ be as in $\ref{A5}$ and let $H$ be an object of $\cA_X$. Then we have a unique direct sum decomposition $H=\bigoplus_c I(c)$ in $\cA$, where $c$ ranges over all elements of $\bar \Q_{\ell}^\times$, such that on the stalk of $I(c)$, 
the action of $\gamma_1-c$ is nilpotent.

\end{lem} 
\begin{pf} Let $V_c$ be the generalized eigenspace of $\gamma_1$ for the eigenvalue $c$ in the stalk of $H$. We prove that $V_c$ is preserved by the action of $\pi_1^{\log}(X)$. The Frobenius $F$ and $\gamma_2$ commute with $\gamma_1$, so they preserve $V_c$. We prove $\gamma_0$ preserves $V_c$. From $\gamma_0\gamma_1\gamma_0^{-1}= \gamma_1\gamma_2$, we have 
 $(\gamma_1-c)^n\gamma_0^{-1}=\gamma_0^{-1}((\gamma_1-c)\gamma_2+c(\gamma_2-1))^n$ for $n\geq 0$. If $n$ is sufficiently large, since the action of $\gamma_2-1$ is nilpotent (Lemma \ref{gamma2}), 
  $((\gamma_1-c)\gamma_2+c(\gamma_2-1))^n$ kills $V_c$ and hence $(\gamma_1-c)^n$  kills $\gamma_0^{-1}V_c$.  Hence $\gamma_0^{-1}$ preserves $V_c$.

The projections to the generalized eigenspaces are given by a polynomial of $\gamma_1$ and hence preserve relative monodromy filtrations, and hence they give a direct decomposition in $\cA_X$. 
\end{pf}

\begin{lem}\label{neq1} Let $X\to Y$ be as in $\ref{A5}$, let $H$ be an object of $\cA_X$, and assume that $H=I(c)$ (Lemma $\ref{eigen}$) with $c\neq 1$. Then $H^i(H)=0$ for all $i$. 

\end{lem}

\begin{pf} As a $\bar \Q_{\ell}$-vector space, the stalk of $H^i(H)$ is identified with the $i$-th cohomology of the complex $0\to H\to H\oplus H \to H\to 0$, where $H$ denotes the stalk of $H$, $H\to H\oplus H$ is $(\gamma_1-1, \gamma_2-1)$, and $H\oplus H\to H$ is $(1-\gamma_2, \gamma_1-1)$. It is acyclic because $\gamma_1-1:H\to H$ is an isomorphism. \end{pf}

\begin{para} We prove Proposition \ref{Hi2}. 

(1) follows from Lemma \ref{eigen}, and (2) follows from Lemma \ref{neq1}.

 We prove (3). 

The implication (i) $\Rightarrow$ (ii) and its version for $\cA^s_Y$ are clear.

 We consider the actions of $N_j=\log(\gamma_j)$ ($j=1,2$) on the stalk of $H$. Note that $N_1$ is of Frobenius weight $0$ and $N_2$ is of Frobenius weight $-2$. 
 
 We prove the implication (ii) $\Rightarrow$ (iv). We have $H^0(H)\neq 0$ because there is a non-zero element $x$ of the lowest Frobenius weight part of the stalk of $H$ such that $N_1(x)=0$ ($N_2(x)=0$ automatically). Since $H^0(H)\in \cA_Y$, we have $H^0(H)\in \cA_X$. Since $H$ is simple, we have $H=H^0(H)$, that is,  (iv) is satisfied. 
 
  We prove the implication (iii) $\Rightarrow$ (iv). As $\bar \Q_{\ell}$-vector spaces, we can identify  $H^0(H)$ with the part $H^{N_1=N_2=0}$ of $H$ (this $H$ denotes the stalk of $H$), $H^2(H)$  with $H/(N_1H+N_2H)$, and the Lefschetz class map with the map induced by the identity map of $H$. Then the assumption that it is an isomorphism 
 implies $N_1=N_2=0$. In fact, if $N_j\neq 0$ for some $j$, $H^{N_j=0} \cap N_jH$ is a non-zero subspace, on which $N_{3-j}$ acts.  Hence $H^{N_1=N_2=0} \cap N_jH\not=\{0\}$. Hence $H$ comes from $Y$. This object on $Y$ belongs to $\cA_Y$ 
  because its pullback $H$  belongs to $\cA_X$.

 The implication (iv) $\Rightarrow$ (i) follows from Proposition  \ref{Hi}. 
 
 We prove the implication (iv) $\Rightarrow$ (iii). 
 This is reduced to the case where $H$ is the constant object
  $\bar \Q_{\ell}$. In this case, $H^0(H)=\bar \Q_{\ell}$ of pure $W$-weight $0$, $H^2(H)= \bar \Q_{\ell}(-1)$ of pure $W$-weight $2$, and the Lefschetz map is identified with the identity map $\bar \Q_{\ell}\to \bar \Q_{\ell}$. 
 \end{para}

\begin{para}

The authors expect that the equivalent  conditions in Proposition \ref{Hi2} (3) are  actually always satisfied. If this is the case, the statement in Question \ref{Hope} is true for $X\to Y$ in \ref{A5}. To see this, by Lemma \ref{eigen}, we may assume $H=I(c)$ for some $c$. If $c\neq 1$, it follows from Lemma \ref{neq1}. If $c=1$, using (iv) in Proposition \ref{Hi2} (3), it follows from Proposition  \ref{Hi}. 

\end{para}

\begin{para}\label{A4} An example for non-semi-simple case which does not have the property in the statement in Question \ref{Hope}.

Let $X\to Y$ be as in \ref{A5}. 
Let $H$ be the two dimensional $\bar \Q_{\ell}$-space with base $(f_j)_j$ endowed with the following action of $\pi_1^{\log}(X)$ over $\bar \Q_{\ell}$. Take a splitting of the exact sequence \ref{pi1}, let the action of $\pi_1^{\log}(y)$ be the trivial one, and let the action of $\pi_1^{\log}(X \times_Y \bar y(\log))$ be the one induced from the action of $T$, where $\gamma_2$ acts trivially and $\gamma_1$ acts as $\gamma_1(f_1)=f_1$, $\gamma_1(f_2)=f_2+f_1$. Let $W$ of $H$ be pure of weight $0$.

  This $H$ is a $W$-pure object of $\cA_X$ but is not semi-simple ($f_1$ generates a subobject in $\cA_X$ which is not a direct summand).  This  $H$ does not satisfy the statement in Question  \ref{Hope}. The first higher direct image $H^1(H)$ of $H$ does not belong to $\cA_Y$.  This is seen as follows.

  Note that  $H^1(H)$ has the pure $W$-weight $1$ by the definition of the $W$ of the higher direct image. We identify $L^*$ with $S_1$ (\ref{loget}). The map $H \to H \otimes L^*$ in \ref{cpx} sends $f_2$ to $f_1\otimes e_1$ and sends $f_1$ to $0$. From this, we see that  $H^1(H)$ 
    is identified with the subquotient $P/Q$  of $H\otimes L^*$, where $P$ has the base $f_2\otimes e_1$, $f_1\otimes e_2$,  $f_1\otimes e_1$ and $Q$ has the base $f_1\otimes e_1$. Hence $H^1(H)$ has the base $f_2\otimes e_1\bmod Q$, $ f_1\otimes e_2\bmod Q$, the former has Frobenius weight $0$ and the latter has Frobenius weight $2$.  Since $\gamma_0(f_j)=f_j$, $\gamma_0(e_1)=e_1$, and $\gamma_0(e_2) \equiv e_2 \bmod \bar \Q_{\ell}e_1$, 
   we have $\gamma_0(f_2 \otimes e_1)=f_2 \otimes e_1$ and  $\gamma_0(f_1\otimes e_2)\equiv f_1\otimes e_2 \bmod Q$, and hence  $\gamma_0$ acts trivially  on $H^1(H)$. Hence the relative monodromy filtration of a non-zero element of $\sig(y)$ (which acts as zero) on $H^1(H)$ is pure of weight $1$, and does not coincide with the Frobenius weight filtration.
  Hence $H^1(H)$ is not an object of $\cA_Y$.

\end{para}

\noindent {\rm Kazuya KATO
\\
Department of mathematics
\\
University of Chicago
\\
Chicago, Illinois, 60637, USA}
\\
{\tt kkato@uchicago.edu}

\bigskip

\noindent {\rm Chikara NAKAYAMA
\\
Department of Economics 
\\
Hitotsubashi University 
\\
2-1 Naka, Kunitachi, Tokyo 186-8601, Japan}
\\
{\tt c.nakayama@r.hit-u.ac.jp}

\bigskip

\noindent
{\rm Sampei USUI
\\
Graduate School of Science
\\
Osaka University
\\
Toyonaka, Osaka, 560-0043, Japan}
\\
{\tt usui@math.sci.osaka-u.ac.jp}

\end{document}